\documentclass[12pt,a4paper]{article}

\usepackage{graphicx}
\usepackage{subfig}
\usepackage{amsmath,amsthm,amssymb}
\usepackage{esint}
\usepackage{mathrsfs}

\usepackage[T1]{fontenc}

\usepackage{hyperref}

\usepackage[left=2.2cm,right=2.2cm,top=3cm,bottom=3.5cm]{geometry}

\usepackage{color}
\usepackage{array}

\usepackage{multirow}

\theoremstyle{plain}
\newtheorem{thm}{Theorem}[section]
\newtheorem{prop}[thm]{Proposition}

\theoremstyle{definition}
\newtheorem{defn}[thm]{Definition}

\theoremstyle{remark}
\newtheorem{rem}[thm]{Remark}

\numberwithin{equation}{section}

\usepackage{csquotes}
\usepackage[backend=bibtex8, style=numeric-comp, maxnames=5, sorting=nyt]{biblatex}
\usepackage{enumitem}
\renewcommand\labelenumi{(\alph{enumi})}
\renewcommand\theenumi\labelenumi

\newcommand{\R}{\mathbb{R}} 

\newcommand{\Grad}{\nabla}  
\newcommand{\Div}{{\rm div}\,} 

\newcommand{\dx}{\,{\rm d}\vx}
\newcommand{\dxc}{\,{\rm d}x}
\newcommand{\dyc}{\,{\rm d}y}
\newcommand{\dt}{\,{\rm d}t}
\newcommand{\dtau}{\,{\rm d}\tau}
\newcommand{\dr}{\,{\rm d}r} 

\renewcommand{\vec}[1]{{\bf #1}}
\newcommand{\vx}{\vec{x}}
\newcommand{\vu}{\vec{u}} 
\newcommand{\vphi}{\boldsymbol{\varphi}} 

\renewcommand{\rho}{\varrho}
\newcommand{\half}{\tfrac{1}{2}}  
\newcommand{\Cc}{C^\infty_{\rm c}} 
\newcommand{\ov}[1]{\overline{#1}}  

\newcommand{\ci}{{\rm ci}}

\newcommand{\tc}{{\rm c}}

\newcommand{\name}[1]{\textsc{#1}} 

\begin{document} 

\renewcommand{\arraystretch}{1.2}

\title{The local least action criterion fails as a selection criterion for weak solutions of the compressible Euler equations} 

\author{Simon Markfelder\footnote{e-mail: \textsf{simon.markfelder@uni-konstanz.de}} \and \begingroup \renewcommand{\thefootnote}{\ensuremath{\dagger,\ddagger}} Valentin Pellhammer \footnote{Corresponding author, e-mail: \textsf{valentin.pellhammer@uni-hohenheim.de}} \endgroup } 

\date{\today}

\maketitle

\bigskip

\centerline{$^\ast$ Universit\"at Konstanz, Department of Mathematics and Statistics,} 
\centerline{Post office box: 199, 78457 Konstanz, Germany} 

\bigskip

\centerline{$^\dagger$ University of Hohenheim, Institute of Applied Mathematics and Statistics,} 
\centerline{70599 Stuttgart, Germany} 

\bigskip

\centerline{$^\ddagger$ Computational Science Hub, University of Hohenheim} 
\centerline{70599 Stuttgart, Germany} 

\bigskip

\begin{abstract} 
	It is now well known that the classical notion of admissible weak solutions (also known as weak entropy solutions) does not restore uniqueness for 
	the multi-dimensional compressible Euler equations. Indeed, convex integration has shown that admissible weak solutions are in general highly non-unique. 
	This has motivated additional selection criteria intended to rule out the counterintuitive solutions generated by convex integration. In this paper we 
	prove that the local least action criterion introduced by H.~Gimperlein, M.~Grinfeld, R.~J.~Knops and M.~Slemrod does not serve as a proper selection 
	criterion, since it fails to select the solution that is intuitively expected to be physically relevant. 
\end{abstract}

\bigskip


\bigskip



\tableofcontents

\section{Introduction} \label{sec:intro} 

This paper deals with the issue of finding a proper notion of solution for the isentropic Euler equations
\begin{align}
	\partial_t \rho + \Div (\rho\vu) &= 0, \label{eq:euler-mass} \\
	\partial_t (\rho\vu) + \Div(\rho\vu\otimes \vu) + \Grad p(\rho) &= 0 \label{eq:euler-mom}
\end{align}
in two and higher space dimensions. To keep the presentation simple, we restrict to two dimensions, but our main result can be trivially generalized to any dimension larger or equal 2. So, the unknowns in \eqref{eq:euler-mass}, \eqref{eq:euler-mom} are the scalar-valued density $\rho$ and the vector-valued velocity $\vu$, both of which are functions of time $t\in [0,\infty)$ and space $\vx= (x,y)\in \R^2$. The solutions considered in this paper will not contain vacuum, i.e.~there will hold $\rho(t,\vx)>0$. The pressure $p$ in \eqref{eq:euler-mom} is a given function of the density $\rho$. For simplicity, we consider the polytropic gas law with adiabatic coefficient $\gamma=2$, i.e.
\begin{equation} \label{eq:pressure}
	p(\rho) = \rho^2. 
\end{equation} 

We are interested in the initial value problem for \eqref{eq:euler-mass}, \eqref{eq:euler-mom}, where the initial state is given by $(\rho_0,\vu_0)\in L^\infty(\R^2;\R^+\times \R^2)$. In other words we endow the Euler equations \eqref{eq:euler-mass}, \eqref{eq:euler-mom} with the initial condition
\begin{equation} \label{eq:init}
	(\rho,\vu) (0,\cdot) = (\rho_0,\vu_0). 
\end{equation}

It is well-known that strong solutions to \eqref{eq:euler-mass}-\eqref{eq:init} do not exist globally in time as shocks may evolve. For this reason, one has to consider weak solutions. However, weak solutions are highly non-unique and allow for the unphysical phenomenon that the energy increases in time. In order to rule out such behavior, one usually studies so-called \emph{admissible weak solutions}, i.e.~weak solutions which in addition satisfy the energy inequality 
\begin{equation} \label{eq:euler-en}
	\partial_t \Big( \half \rho |\vu|^2 + P(\rho) \Big) + \Div\Big[ \Big( \half \rho |\vu|^2 + P(\rho) + p(\rho) \Big) \vu \Big] \leq 0 
\end{equation}
in the sense of distributions. In \eqref{eq:euler-en}, the pressure potential $P$ is given by 
$$
	P(\rho) = \rho \int^\rho \frac{p(r)}{r^2}\dr,
$$
which can be simplified for our choice \eqref{eq:pressure} to $P(\rho)=\rho^2$. For the sake of completeness, we now provide a precise definition of admissible weak solutions. 

\begin{defn}
	A pair $(\rho,\vu)\in L^\infty ((0,\infty)\times\R^2;\R^+\times\R^2)$ is called \emph{weak solution} of the initial value problem \eqref{eq:euler-mass}-\eqref{eq:init} if 
	\begin{align*}
		\int_0^\infty \int_{\R^2}\Big[ \rho \partial_t \phi +\rho \vu \cdot\nabla \phi \Big] \dx \dt + \int_{\R^2}\rho_0\phi(0,\cdot) \dx &=0, \\
		\int_0^\infty \int_{\R^2}\Big[ \rho\vu\cdot\partial_t \vphi +\rho\vu\otimes\vu:\nabla \vphi +p(\rho)\Div \vphi \Big] \dx \dt + \int_{\R^2} \rho_0\vu_0\cdot \vphi(0,\cdot) \dx &=0
	\end{align*}
	hold for all test functions $(\phi,\vphi)\in \Cc([0,\infty)\times \R^2;\R\times \R^2)$. A weak solution $(\rho,\vu)$ is called \emph{admissible} if in addition
	\begin{align*}
		\int_0^\infty \int_{\R^2}\left[ \Big(\half\rho|\vu|^2 + P(\rho)\Big)\partial_t\psi + \Big(\half\rho|\vu|^2 +P(\rho) + p(\rho)\Big)\vu\cdot\nabla\psi \right]\dx\dt & \\
		+\int_{\R^2}\Big(\half\rho_0|\vu_0|^2+P(\rho_0)\Big)\psi(0,\cdot)\dx &\geq 0 
	\end{align*}
	for all $\psi\in\Cc([0,\infty)\times \R^2;\R_0^+)$.
\end{defn}

It is nowadays well-known that for many initial data, the initial value problem \eqref{eq:euler-mass}-\eqref{eq:init} admits infinitely many admissible weak solutions. In other words, the energy inequality \eqref{eq:euler-en} is not able to single out a unique solution. This fact is revealed by the technique of convex integration. Non-uniqueness of admissible weak solutions for the isentropic Euler equations \eqref{eq:euler-mass}, \eqref{eq:euler-mom} was shown for the first time by \name{De Lellis}-\name{Sz{\'e}kelyhidi}~\cite{DelSze10}, a result which is built upon an analogous non-uniqueness statement for the incompressible Euler equations, see also \cite{DelSze09}. Non-uniqueness even occurs for a large class of initial data, see \name{Chiodaroli}~\cite{Chiodaroli14}, \name{Feireisl}~\cite{Feireisl14}, \name{Chiodaroli}-\name{Feireisl}~\cite{ChiFei24_1} and \name{Boutros}-\name{Markfelder}~\cite{BouMar26pre}.

\name{Chiodaroli}-\name{De~Lellis}-\name{Kreml}~\cite{ChiDelKre15} studied a special type of initial data, namely so-called \emph{Riemann initial data} of the form
\begin{equation} \label{eq:riemann}
	(\rho_0,\vu_0)(\vx) = \left\{ \begin{array}{cc} (\rho_-,\vu_-), & \text{ if } y<0, \\ (\rho_+,\vu_+), & \text{ if } y>0. \end{array} \right. 
\end{equation}
For such data one can explicitly compute an admissible weak solution by solving the corresponding one-dimensional Riemann problem in $y$-direction, see e.g.~\cite[Sect.~7.1]{Markfelder} or the standard monograph by \name{Dafermos}~\cite{Dafermos}. We will call this solution the \emph{1-D solution} to the problem \eqref{eq:euler-mass}-\eqref{eq:init}, \eqref{eq:riemann}. As shown in \cite{ChiDelKre15}, there exist initial data of the form \eqref{eq:riemann} such that the initial value problem \eqref{eq:euler-mass}-\eqref{eq:init} has infinitely many additional admissible weak solutions. It turned out that non-uniqueness occurs as soon as the 1-D solution contains a shock, see the series of papers \cite{ChiDelKre15,ChiKre14,ChiKre18,KliMar18_1,BreChiKre18} and also \cite{Markfelder} for a summary, while the 1-D solution is the unique solution if it only contains rarefaction waves, see \name{Chen}-\name{Chen}~\cite{CheChe07} and also \name{Feireisl}-\name{Kreml}~\cite{FeiKre15}. Non-uniqueness in the case where the 1-D solution only consists of a contact discontinuity has been shown recently by \name{Krupa}-\name{Sz{\'e}kelyhidi}~\cite{KruSze25} for special (unphysical) pressure laws, and by \name{Horimoto}~\cite{Horimoto26pre} in the case of the polytropic pressure, in particular for \eqref{eq:pressure}. 

In order to cure the apparent non-uniqueness problem, several selection criteria have been studied in the literature. Those criteria are supposed to rule out unphysical solutions with the hope that only one solution -- the physically relevant solution -- remains. As a test case, one has considered Riemann data of the form \eqref{eq:riemann} which give rise to infinitely many admissible weak solutions. In particular, one has checked if the criterion at hand selects the 1-D solution which, intuitively, should be the physically relevant solution in this context. \name{Chiodaroli}-\name{Kreml}~\cite{ChiKre14} considered the \emph{entropy rate criterion} which picks the solution that dissipates energy at the highest rate, see also \name{Dafermos}~\cite{Dafermos73} who studied the same criterion in a slightly different context. The result, however, is negative: the entropy rate criterion does not select the 1-D solution, which is shown in \cite{ChiKre14} by constructing another solution that dissipates energy at higher rate. A similar statement still holds for a localized version of the entropy rate criterion as proved by the first author in~\cite{Markfelder24}. More recently, the \emph{least action criterion} was proposed in the context of the Euler equations \eqref{eq:euler-mass}, \eqref{eq:euler-mom} by \name{Gimperlein~et~al.}~\cite{GGKS25}. The least action criterion suggests to select the admissible weak solution which minimizes the action functional
$$
	\mathcal{A}[\rho,\vu] := \int_0^T \int_\Omega \Big( \half \rho |\vu|^2 - P(\rho) \Big) \dx\dtau, 
$$
where $T>0$ denotes the final time, and $\Omega$ is the spatial domain\footnote{In the context of Riemann initial data \eqref{eq:riemann}, where $\Omega=\R^2$, one may restrict to a sufficiently large square $[-L,L]^2$ when computing the spatial integral in the action $\mathcal{A}$, see below.}. However, in~\cite{MarPel25pre} the authors constructed a solution with smaller action than the 1-D solution, which shows that the 1-D solution is not singled out by the least action criterion. The example constructed in \cite{MarPel25pre} also reveals that the solution which is preferred depends on the final time $T>0$ up to which the action functional $\mathcal{A}$ is computed. This fact may be interpreted as further evidence that the least action criterion fails to serve as a proper selection criterion. 

As a response to \cite{MarPel25pre}, \name{Gimperlein~et~al.}~\cite{GGKS26} proposed a \emph{local-in-time} version of the least action criterion, which we call the \emph{local least action criterion}, in order to overcome the aforementioned issues. To be able to give a precise formulation of the criterion, we denote the action computed on $[-L,L]^2$ and on the time interval $[0,t]$ for any $t>0$ by
$$
	\mathcal{A}_L[\rho,\vu](t) := \int_0^t \int_{[-L,L]^2} \Big( \half \rho |\vu|^2 - P(\rho) \Big) \dx\dtau. 
$$

\begin{defn}[{cf.~\cite[Defn.~2]{GGKS26}}] \label{defn:LLA}
	An admissible weak solution $(\rho,\vu)$ to \eqref{eq:euler-mass}-\eqref{eq:init}, \eqref{eq:riemann} satisfies the \emph{local least action criterion}\footnote{In \cite{GGKS26} the local least action criterion is called ``LAAP$_0$-criterion''.} if for all other admissible weak solutions $(\widetilde{\rho},\widetilde{\vu})\neq (\rho,\vu)$ there exists a time $t_1>0$ and an $\ov{L}>0$ such that $\mathcal{A}_L[\rho,\vu](t) \leq \mathcal{A}_L[\widetilde{\rho},\widetilde{\vu}](t)$ for all $t\in (0,t_1)$ and all $L\geq \ov{L}$.
\end{defn}

\begin{rem}
	In \cite{GGKS26} the local least action criterion is stated in a more general way, namely on a time interval $(t_0,t_1)$. For our purposes it suffices to take $t_0=0$. 
\end{rem}

\begin{rem}
	The restriction to the square $[-L,L]^2$ is one possibility to deal with the fact the domain $\R^2$ is unbounded. Another option is to periodize as carried out in \cite[Sect.~6]{ChiKre14}.
\end{rem}

\name{Gimperlein~et~al.}~\cite{GGKS26} proved that a large class of solutions constructed by convex integration (including the solution produced in \cite{MarPel25pre}) does not beat the 1-D solution with respect to the order relation induced by Defn.~\ref{defn:LLA}. This suggests that the local least action criterion may indeed be appropriate to select the ``right'' solution. The main purpose of this paper, however, is to give a negative answer. Our main result reads as follows.

\begin{thm} \label{thm:main}
	There exist initial states $(\rho_\pm,\vu_\pm)\in \R^+ \times \R^2$ such that the 1-D solution to \eqref{eq:euler-mass}-\eqref{eq:init}, \eqref{eq:riemann} does not satisfy the local least action criterion. 
\end{thm}

\begin{rem} 
	As a matter of fact, our proof of Thm.~\ref{thm:main} will also show that the 1-D solution does not satisfy the \emph{strict} local least action criterion, and neither the \emph{strict action rate criterion} as defined in \cite[Defns.~2 and 4]{GGKS26}.
\end{rem}

Sect.~\ref{sec:proof} is devoted to the proof of Thm.~\ref{thm:main}. To this end, we will consider an example of initial states $(\rho_\pm,\vu_\pm)\in \R^+ \times \R^2$ and construct a solution which beats the 1-D solution with respect to the order relation induced by Defn.~\ref{defn:LLA}. We finish this paper with a short discussion in Sect.~\ref{sec:discussion}.

\section{Proof of the main result} \label{sec:proof}

In this section we will prove the main result Thm.~\ref{thm:main}. In Sect.~\ref{subsec:pf-ci} we recall a convex integration result from \cite{ChiDelKre15}. A particular pair of initial states $(\rho_\pm,\vu_\pm)\in \R^+ \times \R^2$ is introduced in Sect.~\ref{subsec:pf-states}. Sect.~\ref{subsec:pf-1D} is devoted to the 1-D solution corresponding to $(\rho_\pm,\vu_\pm)$. In Sect.~\ref{subsec:pf-wild} we present another solution and compute its action. We conclude in Sect.~\ref{subsec:pf-conclusion} by comparing the action of the two solutions.

\subsection{Infinitely many solutions for Riemann initial data} \label{subsec:pf-ci}

The following statement was originally shown by \name{Chiodaroli}-\name{De Lellis}-\name{Kreml}~\cite{ChiDelKre15}. 

\begin{prop}\label{prop:ci} 
	Let $\rho_\pm\in\R^+$, $\vu_\pm =(u_\pm,v_\pm)\in\R^2$. Assume that there exist numbers $\mu_0,\mu_1\in\R$, $\rho_1\in\R^+$, $\vu_1=(u_1,v_1)\in \R^2$, $\gamma_1,\delta_1\in\R$ and $C_1\in\R^+$ which fulfill the following algebraic equations and inequalities:
	\begin{itemize} 
		\item Order of speeds:
		$$
			\mu_0<\mu_1;
		$$
		
		\item Rankine Hugoniot conditions on the left interface:
		\begin{align*}
			\mu_0(\rho_- - \rho_1) &= \rho_-v_- - \rho_1v_1 ; \\
			\mu_0(\rho_-u_- - \rho_1u_1) &= \rho_-u_-v_- - \rho_1\delta_1 ; \\
			\mu_0(\rho_-v_--\rho_1v_1) &=\rho_-(v_-)^2 - \rho_1\left(\frac{C_1}{2}-\gamma_1\right)+p( \rho_-) - p(\rho_1);
		\end{align*}
		
		\item Rankine Hugoniot conditions on the right interface:
		\begin{align*}
			\mu_1(\rho_1 - \rho_+) &= \rho_1v_1 - \rho_+v_+ ;\\
			\mu_1(\rho_1u_1 - \rho_+u_+) &= \rho_1\delta_1 - \rho_+u_+v_+ ;\\
			\mu_1(\rho_1v_1-\rho_+v_+) &=\rho_1\left(\frac{C_1}{2}-\gamma_1\right) - \rho_+(v_+)^2 + p(\rho_1) - p(\rho_+);
		\end{align*}
		
		\item Subsolution condition:
		\begin{align*}
			C_1 - (u_1)^2 - (v_1)^2&>0; \\
			\left(\frac{C_1}{2} - (u_1)^2 + \gamma_1\right)\left(\frac{C_1}{2} - (v_1)^2 - \gamma_1\right) - (\delta_1 - u_1v_1)^2&>0 \notag ;
		\end{align*}
		
		\item Admissibility condition on the left interface:
		\begin{align*}
			&\mu_0\left(\rho_-\frac{|\vu_-|^2}{2} + P(\rho_-) - \rho_1\frac{C_1}{2}-P(\rho_1)\right)\\
			&\leq \left(\rho_-\frac{|\vu_-|^2}{2} + P(\rho_-) + p(\rho_-)\right)v_- - \left(\rho_1\frac{C_1}{2} + P(\rho_1) + p(\rho_1)\right)v_1;
		\end{align*}
		
		\item Admissibility condition on the right interface:
		\begin{align*}
			&\mu_1\left(\rho_1 \frac{C_1}{2} + P(\rho_1) - \rho_+\frac{|\vu_+|^2}{2} - P(\rho_+)\right)\\
			&\leq \left(\rho_1\frac{C_1}{2}+P(\rho_1)+p(\rho_1)\right)v_1 - \left(\rho_+\frac{|\vu_+|^2}{2} + P(\rho_+) + p(\rho_+)\right)v_+.
		\end{align*}
	\end{itemize}		
	Then there exist infinitely many admissible weak solutions $(\rho,\vu)\in L^\infty((0,\infty)\times \R^2; \R^+ \times \R^2)$ of the initial value problem \eqref{eq:euler-mass}-\eqref{eq:init}, \eqref{eq:riemann}. All those solutions fulfill
	\begin{itemize}
		\item $\rho(t,\vx)= \rho_1$ and $|\vu(t,\vx)|^2= C_1$ for a.e.~$(t,\vx)\in \Gamma_1$, 
		\item $(\rho,\vu)(t,\vx) = (\rho_\pm,\vu_\pm)$ for a.e.~$(t,\vx)\in \Gamma_\pm$, 
	\end{itemize}
	where $\Gamma_-,\Gamma_1,\Gamma_+$ are given by 
	\begin{align*}
		\Gamma_- &:= \left\{ (t,\vx)\in (0,\infty)\times\R^2\,\big|\, y < \mu_0 t \right\}, \\
		\Gamma_1 \,&:= \left\{ (t,\vx)\in (0,\infty)\times\R^2\,\big|\, \mu_0 t < y < \mu_1 t \right\}, \\
		\Gamma_+ &:= \left\{ (t,\vx)\in (0,\infty)\times\R^2\,\big|\, y > \mu_1 t \right\}.
	\end{align*}
\end{prop} 

For the proof of Prop.~\ref{prop:ci} we refer to \cite[Props.~3.6 and 5.1]{ChiDelKre15}, \cite[Props.~3.1 and 4.1]{ChiKre14} or \cite[Thm.~7.3.4 and Prop.~7.3.5]{Markfelder}.

\subsection{A particular pair of initial states} \label{subsec:pf-states}

In order to prove Thm.~\ref{thm:main}, we consider the following initial states:
\begin{align} \label{eq:initial-states}
	\rho_- &= 1, & \vu_- &= \left(\begin{array}{c} -6 \\ \frac{12\sqrt{14}}{7} \end{array}\right), & \rho_+ &= 2, & \vu_+ &= \left(\begin{array}{c} 6 \\ -\frac{15\sqrt{14}}{14} \end{array}\right). 
\end{align}

\subsection{The 1-D solution and its action} \label{subsec:pf-1D} 

Using e.g.~\cite[Prop.~7.1.1]{Markfelder}, it is not difficult to verify that the 1-D solution $(\rho_{1d},\vu_{1d})$ which corresponds to the initial states \eqref{eq:initial-states} consists of a 1-shock, a 2-contact discontinuity, and a 3-shock. The intermediate states are given by 
\begin{align} \label{eq:1-D-sol_middle-states} 
	\rho_{M} &= 7, & \vu_M^- &= \left(\begin{array}{c} -6 \\ 0 \end{array}\right), & \vu_M^+ &= \left(\begin{array}{c} 6\\0 \end{array}\right), 
\end{align}
while the speeds of the shock waves read
\begin{align} \label{eq:1-D-sol_shock-speeds}
	\sigma_- &= -\frac{2\sqrt{14}}{7}, & \sigma_+ &= \frac{3\sqrt{14}}{7}, 
\end{align}
and the speed of the contact discontinuity coincides with the second component of the intermediate velocities, i.e.~$\sigma_\tc = 0$. The 1-D solution $(\rho_{1d},\vu_{1d})$ is sketched in Fig.~\ref{fig:solution-1d} below. 

Let us next compute the action of the 1-D solution $(\rho_{1d},\vu_{1d})$. For $L>\max\{-\sigma_-,\sigma_+\}$ we obtain 
\begin{align*}
	\mathcal{A}_L[\rho_{1d},\vu_{1d}](t) &= \int_0^t \int_{-L}^L \int_{-L}^L \Big( \half \rho_{1d} |\vu_{1d}|^2 - P(\rho_{1d}) \Big) \dxc\dyc\dt \\
	&= 2L \int_0^t \left( \int_{-L}^{\sigma_- \tau} a_-  \dyc + \int_{\sigma_- \tau}^{0} a_M^- \dyc + \int_{0}^{\sigma_+ \tau} a_M^+ \dyc + \int_{\sigma_+ \tau}^L a_+ \dyc \right) \dtau,
\end{align*}
where, according to \eqref{eq:1-D-sol_middle-states},
\begin{align} \label{eq:1-D-sol_action}
	a_- &= \frac{263}{7}, & a_M^- &=77, &
	a_M^+ &=77, & a_+ &=\frac{673}{14}. 
\end{align}
Inserting the latter values as well as \eqref{eq:1-D-sol_shock-speeds}, we end up with
\begin{equation} \label{eq:action-1d}
	\mathcal{A}_L[\rho_{1d},\vu_{1d}](t) = L \left( \frac{2319 \sqrt{14}}{98} t^2 + \frac{1199}{7} L t \right).
\end{equation}

\subsection{Another solution and its action} \label{subsec:pf-wild} 

Next, we invoke Prop.~\ref{prop:ci} to produce another solution $(\rho_\ci,\vu_\ci)$ to \eqref{eq:euler-mass}-\eqref{eq:init}, \eqref{eq:riemann} with initial states \eqref{eq:initial-states}. 

\begin{rem} 
	Note that initial states which give rise to a 1-D solution consisting of two shocks and a contact discontinuity, in particular the initial states \eqref{eq:initial-states}, have been studied by \name{B{\v r}ezina}-\name{Chiodaroli}-\name{Kreml}~\cite{BreChiKre18}. It was shown in \cite{BreChiKre18} that such data give rise to infinitely many admissible weak solutions. Mimicking the structure of the 1-D solution, the additional solutions constructed in \cite{BreChiKre18} consist of two ``wild wedges'' $\Gamma_1,\Gamma_2$. We will produce further solutions which only contain one ``wild wedge'' $\Gamma_1$, cf.~Prop.~\ref{prop:ci}. 
\end{rem}

\begin{figure}[h] 
	\centering
	\subfloat[The 1-D solution $(\rho_{1d},\vu_{1d})$.\label{fig:solution-1d}]{ 
		\centering
		\includegraphics[width=0.43\textwidth]{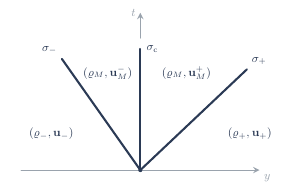}
	}
	\hspace{1.0cm}
	\subfloat[The solution $(\rho_\ci,\vu_\ci)$ constructed by convex integration.\label{fig:solution-ci}]{ 
		\centering 
		\includegraphics[width=0.43\textwidth]{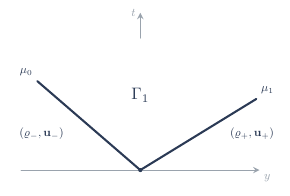} 
	} 
	\caption{Sketch of the solutions $(\rho_{1d},\vu_{1d})$ and $(\rho_\ci,\vu_\ci)$.} 
	\label{fig:solutions}
\end{figure}

We set 
\begin{align*}
	\mu_0 &= -\frac{27\sqrt{14}}{7} + \frac{15\sqrt{35}}{7}, & \mu_1 &= -\frac{27\sqrt{14}}{7} + \frac{20\sqrt{35}}{7}, \\
	\rho_1 &= 5, & C_1 &= \frac{1434052149011}{1914158750} - \frac{43384598136\sqrt{10}}{191415875}, \\
	u_1 &= \frac{66}{5} - \frac{468\sqrt{10}}{125}, & v_1 &= -\frac{96\sqrt{14}}{35} + \frac{12\sqrt{35}}{7}, \\
	\gamma_1 &= \frac{290807163011}{3828317500} - \frac{3972658068\sqrt{10}}{191415875}, & \delta_1 &= -\frac{684\sqrt{14}}{7} + \frac{52272\sqrt{35}}{875}.
\end{align*}
It is then straightforward to check that the algebraic equations and inequalities stated in Prop.~\ref{prop:ci} hold. Thus, Prop.~\ref{prop:ci} yields a solution\footnote{In fact, there are infinitely many such solutions. For our purposes, one of them suffices.} $(\rho_\ci,\vu_\ci)$ to \eqref{eq:euler-mass}-\eqref{eq:init}, \eqref{eq:riemann} with initial states \eqref{eq:initial-states}, which has the properties stated in Prop.~\ref{prop:ci}. The solution $(\rho_\ci,\vu_\ci)$ is sketched in Fig.~\ref{fig:solution-ci}. 

The properties of $(\rho_\ci,\vu_\ci)$ which are mentioned in Prop.~\ref{prop:ci} allow to compute its action. For $L>\max\{-\mu_0,\mu_1\}$ we get 
\begin{align*}
	\mathcal{A}_L[\rho_\ci,\vu_\ci](t) &= \int_0^t \int_{-L}^L \int_{-L}^L \Big( \half \rho_\ci |\vu_\ci|^2 - P(\rho_\ci) \Big) \dxc\dyc\dt \\
	&= 2L \int_0^t \left( \int_{-L}^{\mu_0 \tau} a_-  \dyc + \int_{\mu_0 \tau}^{\mu_1 \tau} a_1 \dyc + \int_{\mu_1 \tau}^L a_+ \dyc \right) \dtau,
\end{align*}
with $a_\pm$ as in \eqref{eq:1-D-sol_action}, and
$$
	a_1 = \frac{1}{2} \rho_1 C_1 - \rho_1^2 = \frac{1414910561511}{765663500} - \frac{21692299068\sqrt{10}}{38283175}. 
$$
Hence, we find using the values for $\mu_0,\mu_1$
\begin{equation} \label{eq:action-ci}
	\mathcal{A}_L[\rho_\ci,\vu_\ci](t) = L \left(\frac{193426517573 \sqrt{35}}{153132700} t^2 -\frac{42516335727 \sqrt{14}}{21438578} t^2 + \frac{1199}{7}L t\right). 
\end{equation}

\subsection{Comparison of the action and conclusion} \label{subsec:pf-conclusion} 

\begin{figure}[h] 
	\centering
	\includegraphics[width=0.7\textwidth]{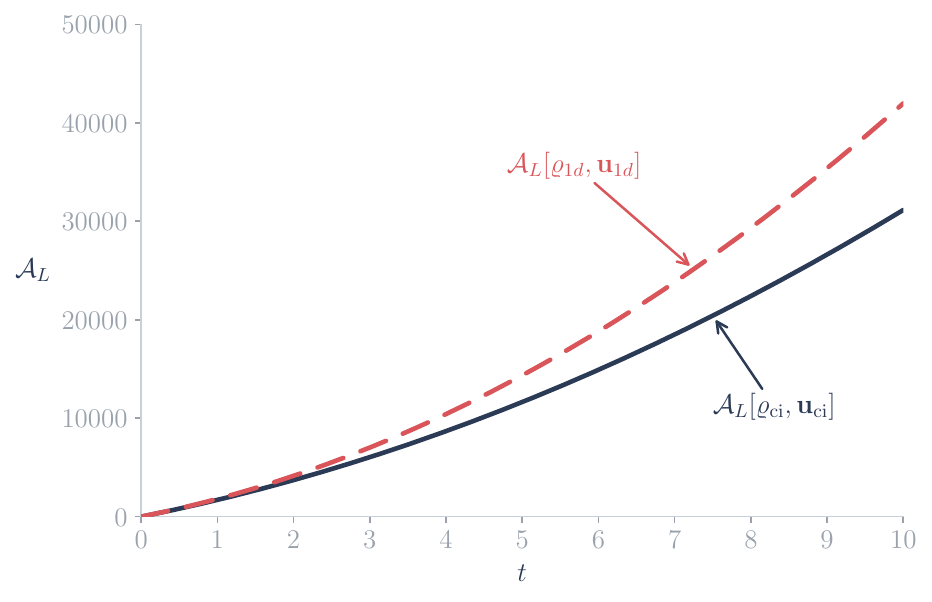} 
	\caption{Comparison of the action $\mathcal{A}_L[\rho_{1d},\vu_{1d}](t)$ and $\mathcal{A}_L[\rho_\ci,\vu_\ci](t)$ for $L=3$.}
	\label{fig:actions}
\end{figure}

For all $L>\max\{-\sigma_-,\sigma_+, -\mu_0,\mu_1\}$ we may compare the action of the 1-D solution $(\rho_{1d},\vu_{1d})$ with the action of the solution $(\rho_\ci,\vu_\ci)$. Using \eqref{eq:action-1d} and \eqref{eq:action-ci}, we can plot the actions, see Fig.~\ref{fig:actions}. This plot (i.e.~Fig.~\ref{fig:actions}) suggests that
\begin{equation} \label{eq:action-comp}
	\mathcal{A}_L[\rho_{1d},\vu_{1d}](t) > \mathcal{A}_L[\rho_\ci,\vu_\ci](t) \qquad \text{ for all } t>0.
\end{equation}

In order to prove \eqref{eq:action-comp} rigorously, we compute 
$$
		\mathcal{A}_L[\rho_{1d},\vu_{1d}](t) - \mathcal{A}_L[\rho_\ci,\vu_\ci](t) = \frac{2151182124300 \sqrt{14}-1353985623011 \sqrt{35} }{1071928900} L t^2,
$$
where we again made use of \eqref{eq:action-1d} and \eqref{eq:action-ci}. As
$$
	 \frac{2151182124300 \sqrt{14}-1353985623011 \sqrt{35} }{1071928900} > 0,
$$
we indeed deduce \eqref{eq:action-comp}, which immediately reveals that the 1-D solution $(\rho_{1d},\vu_{1d})$ does not satisfy the local least action criterion. Thus, the proof of Thm.~\ref{thm:main} is finished.

\section{Discussion} \label{sec:discussion}

As mentioned in Sect.~\ref{sec:intro}, we already know from \cite{ChiKre14}, \cite{Markfelder24} and \cite{MarPel25pre}, respectively, that the \emph{entropy rate criterion}, a localized version of the latter, as well as the \emph{least action criterion} discard the 1-D solution. The main result presented in this paper, i.e.~Thm.~\ref{thm:main}, shows that the \emph{local least action criterion} demonstrates the same behavior, i.e.~it rules out the 1-D solution. 

Note that the example constructed in Sect.~\ref{sec:proof} in order to prove Thm.~\ref{thm:main} represents a further example which shows that the (non-local) least action criterion does not select the 1-D solution. In fact, this example is much simpler than the one produced in \cite{MarPel25pre}. From this point of view, one may consider \cite{MarPel25pre} as obsolete. However, as mentioned in Sect.~\ref{sec:intro}, the example constructed in \cite{MarPel25pre} discloses that the solution which is preferred by the least action criterion depends on the final time $T>0$ up to which the action functional $\mathcal{A}$ is computed. In contrast to that, the example studied in Sect.~\ref{sec:proof} does not exhibit such behavior. 

As a consequence of Thm.~\ref{thm:main}, one either has to reconsider one's intuition that the 1-D solution is the physically relevant solution in the context under consideration, or the local least action criterion must be discarded.

\section*{Acknowledgement}
Funded by the Deutsche Forschungsgemeinschaft (DFG, German Research Foundation) -- SPP 2410 \emph{Hyperbolic Balance Laws in Fluid Mechanics: Complexity, Scales, Randomness (CoScaRa)}, within the Project 525935467 \emph{Convex integration: towards a mathematical understanding of turbulence, Onsager conjectures and admissibility criteria}. 

\printbibliography[heading=bibintoc]

@article{BreChiKre18,
	author  = "J.~B{\v r}ezina and E.~Chiodaroli and O.~Kreml",
	title   = "Contact discontinuities in multi-dimensional isentropic {E}uler equations",
	journal = "Electron. J. Differential Equations",
	year    = "2018",
	volume  = "2018",
	number  = "94",
	pages   = "1-11"
}

@article{CheChe07,
	author  = "G.-Q.~Chen and J.~Chen",
	title   = "Stability of rarefaction waves and vacuum states for the multidimensional {E}uler equations",
	journal = "J. Hyperbolic Differ. Equ.",
	year    = "2007",
	volume  = "4",
	number  = "1",
	pages   = "105-122"
}

@article{Chiodaroli14,
	author  = "E.~Chiodaroli",
	title   = "A counterexample to well-posedness of entropy solutions to the compressible {E}uler system",
	journal = "J. Hyperbolic Differ. Equ.",
	year    = "2014",
	volume  = "11",
	number  = "3",
	pages   = "493-519"
}

@article{ChiDelKre15,
	author  = "E.~Chiodaroli and C.~{De~Lellis} and O.~Kreml",
	title   = "Global ill-posedness of the isentropic system of gas dynamics",
	journal = "Comm. Pure Appl. Math.",
	year    = "2015",
	volume  = "68",
	number  = "7",
	pages   = "1157-1190"
}

@article{ChiFei24_1,
	author  = "E.~Chiodaroli and E.~Feireisl",
	title   = "On the density of ``wild'' initial data for the barotropic {E}uler system",
	journal = "Ann. Mat. Pura Appl. (4)",
	year    = "2024",
	volume  = "203",
	number  = "4",
	pages   = "1809-1817"
}

@article{ChiKre14,
	author  = "E.~Chiodaroli and O.~Kreml",
	title   = "On the energy dissipation rate of solutions to the compressible isentropic {E}uler system",
	journal = "Arch. Ration. Mech. Anal.",
	year    = "2014",
	volume  = "214",
	number  = "3",
	pages   = "1019-1049"
}

@article{ChiKre18,
	author  = "E.~Chiodaroli and O.~Kreml",
	title   = "Non-uniqueness of admissible weak solutions to the {R}iemann problem for isentropic {E}uler equations",
	journal = "Nonlinearity",
	year    = "2018",
	volume  = "31",
	number  = "4",
	pages   = "1441-1460"
}

@article{Dafermos73,
	author  = "C.~Dafermos",
	title   = "The entropy rate admissibility criterion for solutions of hyperbolic conservation laws",
	journal = "J. Differential Equations",
	year    = "1973",
	volume  = "14",
	pages   = "202-212"
}

@book{Dafermos,
	author  = "C.~Dafermos",
	title   = "Hyperbolic conservation laws in continuum physics",
	publisher = "Springer",
	year    = "2016",
	address = "Berlin",
	edition = "4",
	series  = "Grundlehren der mathematischen Wissenschaften",
	number  = "325"
}

@article{DelSze09,
	author  = "C.~{De~Lellis} and L.~{Sz{\'e}kelyhidi~Jr.}",
	title   = "The {E}uler equations as a differential inclusion",
	journal = "Ann. of Math. (2)",
	year    = "2009",
	volume  = "170",
	number  = "3",
	pages   = "1417-1436"
}

@article{DelSze10,
	author  = "C.~{De~Lellis} and L.~{Sz{\'e}kelyhidi~Jr.}",
	title   = "On admissibility criteria for weak solutions of the {E}uler equations",
	journal = "Arch. Ration. Mech. Anal.",
	year    = "2010",
	volume  = "195",
	number  = "1",
	pages   = "225-260"
}

@article{Feireisl14,
	author  = "E.~Feireisl",
	title   = "Maximal dissipation and well-posedness for the compressible {E}uler system",
	journal = "J. Math. Fluid Mech.",
	year    = "2014",
	volume  = "16",
	pages   = "447-461"
}

@article{FeiKre15,
	author  = "E.~Feireisl and O.~Kreml",
	title   = "Uniqueness of rarefaction waves in multidimensional compressible {E}uler system",
	journal = "J. Hyperbolic Differ. Equ.",
	year    = "2015",
	volume  = "12",
	number  = "3",
	pages   = "489-499"
}

@article{GGKS25,
	author  = "H.~Gimperlein and M.~Grinfeld and R.~J.~Knops and M.~Slemrod",
	title   = "The least action admissibility principle",
	journal = "Arch. Ration. Mech. Anal.",
	year    = "2025",
	volume  = "249",
	number  = "2",
	note    = "Paper No. 22"
}

@article{GGKS26,
	author  = "H.~Gimperlein and M.~Grinfeld and R.~J.~Knops and M.~Slemrod",
	title   = "On action rate admissibility criteria",
	journal = "Z. Angew. Math. Phys.",
	year    = "2026",
	volume  = "77",
	number  = "2",
	note    = "Paper No. 57"
}

@article{KliMar18_1,
	author  = "C.~Klingenberg and S.~Markfelder",
	title   = "The {R}iemann problem for the multidimensional isentropic system of gas dynamics is ill-posed if it contains a shock",
	journal = "Arch. Ration. Mech. Anal.",
	year    = "2018",
	volume  = "227",
	number  = "3",
	pages   = "967-994"
}

@article{KruSze25,
	author  = "S.~G.~Krupa and L.~{Sz{\'e}kelyhidi~Jr.}",
	title   = "Contact discontinuities for 2-{D} isentropic {E}uler are unique in 1-{D} but wildly non-unique otherwise",
	journal = "Comm. Math. Phys.",
	year    = "2025",
	volume  = "406",
	number  = "5",
	note    = "Paper No. 109"
}

@article{Markfelder24,
	author  = "S.~Markfelder",
	title   = "A new convex integration approach for the compressible {E}uler equations and failure of the local maximal dissipation criterion",
	journal = "Nonlinearity",
	year    = "2024",
	volume  = "37",
	number  = "11",
	pages   = "1-60"
}

@book{Markfelder,
	author  = "S.~Markfelder",
	title   = "Convex Integration Applied to the Multi-Dimensional Compressible Euler Equations",
	publisher = "Springer",
	year    = "2021",
	address = "Cham, Switzerland",
	series  = "Lecture Notes in Mathematics",
	number  = "2294"
}

@misc{BouMar26pre,
	author  = "D.~W.~Boutros and S.~Markfelder",
	title   = "Non-uniqueness of global-in-time admissible weak solutions to the isentropic compressible {E}uler equations for a dense set of initial data",
	archivePrefix = "arXiv",
	eprint  = "2606.14308",
	year    = "2026"
}

@misc{Horimoto26pre,
	author  = "K.~Horimoto",
	title   = "Non-uniqueness of admissible weak solutions to the two-dimensional barotropic compressible {E}uler system with contact discontinuities",
	archivePrefix = "arXiv",
	eprint  = "2603.23921",
	year    = "2026"
}

@misc{MarPel25pre,
	author  = "S.~Markfelder and V.~Pellhammer",
	title   = "Failure of the least action admissibility principle in the context of the compressible {E}uler equations",
	note  = "To appear in SIAM J. Math. Anal.",
	archivePrefix = "arXiv",
	eprint  = "2502.09292",
	year    = "2026"
}

\end{document}